\documentclass[review,3p,authoryear]{elsarticle}
\usepackage{amssymb}
\usepackage{amsmath}
\usepackage{amsthm}
\usepackage{enumitem}
\usepackage{fancyhdr}
\usepackage{hyperref}
\hypersetup{%
  colorlinks=true,linkcolor=blue,urlcolor=blue,citecolor=blue,%
}
\newcommand{\rtn}{\mathrm{\mathbf{R}}}%-------------------- real number
\newcommand{\N}{\mathrm{\mathbf{N}}}%---------------------- positive integers
\newcommand*{\PR}{\mathrm{\mathbf{P}}}%-------------------- probability
\newcommand*{\EX}{\mathrm{\mathbf{E}}}%-------------------- expectation

\newcommand{\EXlr}[1]{\EX\left[#1\right]}
\newcommand*{\dif}{\,\mathrm{d}}%-------------------------- differential
\newcommand*{\F}{\mathcal{F}}
\newcommand*{\CH}{\mathcal{H}}
\newcommand{\s}{\mathcal{S}}
\newcommand{\BS}{\mathrm{\mathbf{S}}}
\newcommand{\M}{\mathrm{M}}
\newcommand{\h}{\mathbf{H}}
\newcommand*{\prs}{\dif\PR-\mathrm{a.s.}}
\newcommand*{\pts}{\dif\PR\times\dif t-\mathrm{a.e.}}
\newcommand{\Ito}{It\^{o}'s}
\newcommand{\tT}[1][0]{[#1,T]}
\newcommand{\intT}[2][T]{\int^{#1}_{#2}}
\newcommand{\me}{\mathrm{e}}
\newcommand{\one}[1]{{\bf 1}_{#1}}
\DeclareMathOperator*{\sgn}{sgn}
\allowdisplaybreaks[4]
\newtheorem{thm}{Theorem}
\newtheorem{lem}[thm]{Lemma}
\newtheorem{pro}[thm]{Proposition}
\newtheorem{rmk}[thm]{Remark}
\newtheorem{dfn}[thm]{Definition}

\newtheorem{ex}[thm]{Example}
\newproof{pf}{Proof}
\newproof{pou}{{\bf Proof of the uniqueness part}}
\newproof{poe}{{\bf Proof of the existence part}}

\makeatletter
\def\ps@pprintTitle{%
     \let\@oddhead\@empty
     \let\@evenhead\@empty
     \def\@oddfoot{}%
     \let\@evenfoot\@oddfoot}
\makeatother
\begin{document}
\begin{frontmatter}


\title{{\bf Multidimensional BSDEs with uniformly continuous generators and general time
  intervals}
  \tnoteref{found}}
  \tnotetext[found]{Supported by the National Natural Science Foundation of
  China (No. 11371362) and the Fundamental Research Funds for the Central
  Universities (No. 2013RC20).}

\author[cumt,fudan]{Shengjun FAN\corref{cor1}}%
\ead{f\_s\_j@126.com}
\cortext[cor1]{Corresponding author}

\author[cumt]{Lishun XIAO}

\author[cumt]{Yanbin WANG}

\address[cumt]{College of Science, China University of Mining and Technology,
  Xuzhou, Jiangsu, 221116, P.R. China}
\address[fudan]{School of Mathematical Sciences, Fudan University, Shanghai, 20043, P.R. China}

\begin{abstract}
This paper is devoted to solving a multidimensional backward stochastic differential equation
with a general time interval, where the generator
is uniformly continuous in $(y,z)$ non-uniformly with respect to $t$.
By establishing some results on deterministic backward differential equations with general
time intervals, and by virtue of Girsanov's theorem and convolution technique, we establish a new existence
and uniqueness result for solutions of this kind of backward stochastic differential
equations, which extends the results of \citet*{Hamadene2003Bernoulli} and \citet*{FanJiangTian2011SPA}
to the general time interval case.
\end{abstract}

\end{frontmatter}

%\pagestyle{fancy}
%\setlength{\headwidth}{\textwidth}
%\fancyhead{}
%\fancyfoot{}
%\fancyhead[CO]{\footnotesize Multidimensional BSDEs with an infinite time
%  interval and uniformly continuous generators}
%\fancyhead[CE]{\footnotesize Xiao, Fan and Wang}
%\fancyhead[RO,LE]{\footnotesize\thepage}
%\renewcommand{\headrulewidth}{0pt}

\section{Introduction}
\label{sec:Introduction}
In this paper, we are concerned with the following multidimensional backward stochastic
differential equation (BSDE for short in the remaining):
\begin{equation}\label{eq:BSDEs}
  y_t=\xi+\intT{t} g(s,y_s,z_s)\dif s-\intT{t}z_s\dif B_s, \quad t\in\tT,
\end{equation}
where $T$ satisfies $0\leq T\leq +\infty$ called the terminal time;
$\xi$ is a $k$-dimensional random vector called the terminal condition;
the random function
$g(\omega,t,y,z):\Omega\times\tT\times\rtn^k\times\rtn^{k\times d}\mapsto\rtn^k$
is progressively measurable for each $(y,z)$, called the generator of BSDE
\eqref{eq:BSDEs}; and $B$ is a $d$-dimensional Brownian motion. The solution
$(y_t,z_t)_{t\in\tT}$ is a pair of adapted processes. The triple $(\xi,T,g)$ is
called the parameters of BSDE \eqref{eq:BSDEs}. We also denote by BSDE $(\xi,T,g)$
the BSDE with the parameters $(\xi,T,g)$.

The nonlinear  BSDEs were initially introduced by
\citet*{PardouxPeng1990SCL}. They proved an existence and uniqueness result for
solutions of multidimensional BSDEs under the assumptions that the generator $g$ is
Lipschitz continuous in $(y,z)$ uniformly with respect to $t$, where the
terminal time $T$ is a finite constant. Since then, BSDEs have attracted more and more
interesting and many applications on BSDEs have been found in mathematical finance,
stochastic control, partial differential equations and so on (See
\citet*{ElKarouiPengQuenez1997MF} for details).

Many works including \citet*{Mao1995SPA}, \citet*{LepeltierSanMartin1997SPL},
\citet*{Kobylanski2000AP}, \citet*{Bahlali2001CRASPSI}, \citet*{Hamadene2003Bernoulli},
\citet*{BriandLepetierSanMrtin2007Bernoulli}, \citet*{WangHuang2009SPL} and\linebreak[4]
\citet*{FanJiangDavison2010CRASSI},
see also the references therein, have weakened the Lipschitz condition on the generator $g$.
In particular, by virtue of some results on deterministic backward differential
equations (DBDEs for short in the remaining), \citet*{Hamadene2003Bernoulli} proved the existence for solutions
of multidimensional BSDEs when the generator $g$ is uniformly continuous in $(y,z)$.
Furthermore, by establishing an estimate for a linear-growth function,
\citet*{FanJiangDavison2010CRASSI} obtained the uniqueness result under
the same assumptions as those in \citet*{Hamadene2003Bernoulli}. It should be pointed
out that all these works mentioned above only deal the BSDEs with finite time intervals.

\citet*{ChenWang2000JAMSA} first extended the terminal time to the general
case and proved the existence and uniqueness for solutions of
BSDEs under the assumptions that the generator $g$ is Lipschitz
continuous in $(y,z)$ non-uniformly with respect to $t$, which improves
the result of \citet*{PardouxPeng1990SCL} to the infinite time interval case. Furthermore,
\citet*{FanJiang2010SPL} and \citet*{FanJiangTian2011SPA} relaxed the Lipschitz condition
of \citet*{ChenWang2000JAMSA} and obtained two existence and uniqueness
results for solutions of BSDEs with general time intervals, which generalizes
the results of \citet*{Mao1995SPA} and \citet*{LepeltierSanMartin1997SPL} respectively.
Recently, \citet*{HuaJiangShi2013JKSS} extended further the result in \citet*{FanJiang2010SPL}
to the reflected BSDEs case. However, up to now, the question of the existence and uniqueness for solutions of multidimensional BSDEs
with general time intervals and uniformly continuous generators in $(y,z)$
has not been studied.

In this paper, by establishing some results on solutions of DBDEs with general time intervals
and by virtue of Girsanov's theorem and convolution technique, we put forward and prove a general existence and
uniqueness result for solutions
of multidimensional BSDEs with general time intervals and uniformly
continuous generators in $(y,z)$ (see Theorem \ref{thm:MainResult} in Section
\ref{sec:MainResultOnBSDEs}), which extends the results of
\citet*{Hamadene2003Bernoulli} and \citet*{FanJiangDavison2010CRASSI}
to the general time interval case. It should be mentioned that the uniform continuous assumptions for the generator
are not necessarily uniform with respect to $t$ in this result.

We would like to mention that some new troubles arise naturally when we change the terminal time of the BSDE and the DBDE from the finite case to the general case.
For example,
in the case of $T=+\infty$, the integration of a constant over $\tT$ is not
finite any more, $\intT{0}u(t)\dif t\leq C\sup_{t\in\tT}u(t)$ may not hold any longer,
and $\intT{0}v^2(s)\dif s<+\infty$ can not imply $\intT{0}v(s)\dif s<+\infty$.
All these troubles are well overcome in this paper. Furthermore, although the whole idea of the proof
for the existence and uniqueness
of Theorem \ref{thm:MainResult} originates from \citet*{Hamadene2003Bernoulli} and
\citet*{FanJiangDavison2010CRASSI} respectively,
some different arguments from those employed in \citet*{Hamadene2003Bernoulli} is used
to prove the existence part of Theorem \ref{thm:MainResult}. More specifically, in the Step 1 of the proof for the existence part of Theorem \ref{thm:MainResult}, the proof of Lemma \ref{lem:GnSequenceConvolution} is completely different from that of the corresponding result in \citet*{Hamadene2003Bernoulli}, and we
do not use the iteration technique used in \citet*{Hamadene2003Bernoulli} for solutions of BSDE $(\xi,T,g^n)$ (see
\eqref{eq:BSDEsYnInExistenceProof} in Section \ref{sec:ProofOfTheorem}). In addition, the Step 3 of our proof for the
existence part is also very different from that in \citet*{Hamadene2003Bernoulli}.
As a result, the proof procedure is simplified at certain degree.

This paper is organized as follows. Section \ref{sec:Preliminaries}
introduces some usual notations and establishes some results on the solutions of
DBDEs with general time intervals.
Section \ref{sec:MainResultOnBSDEs} is devoted to stating the existence and uniqueness
result on BSDEs --- Theorem \ref{thm:MainResult}. Section \ref{sec:ProofOfTheorem}
gives the detailed proof of Theorem \ref{thm:MainResult}, and Appendix provides the
proof of the results on DBDEs studied in Section \ref{sec:Preliminaries}.

\section{Notations and some results on DBDEs}
\label{sec:Preliminaries}

First of all, let $(\Omega,\F,\PR)$ be a probability space carrying a standard
$d$-dimensional Brownian motion $(B_t)_{t\geq 0}$ and let $(\F_t)_{t\geq 0}$
be the natural $\sigma$-algebra filtration generated by $(B_t)_{t\geq 0}$.
We assume that $\F_T=\F$ and $(\F_t)_{t\geq 0}$ is right-continuous and complete.
In this paper, the Euclidean norm of a vector $y\in \rtn^k$ will be defined by
$|y|$, and for a $k\times d$ matrix $z$, we define $|z|=\sqrt{Tr(zz^*)}$,
where and hereafter $z^*$ represents the transpose of $z$. Let $\langle x,y\rangle$
represent the inner product of $x$, $y\in\rtn^k$.

Let $L^2(\Omega,\F_T,\PR;\rtn^k)$ be the set of $\rtn^k$-valued and $\F_T$-measurable
random variables $\xi$ such that $\|\xi\|^2_{L^2}:=\EX[|\xi|^2]<+\infty$ and let
${\s}^2(0,T;\rtn^k)$ denote the set of
$\rtn^k$-valued, adapted and continuous processes $(Y_t)_{t\in\tT}$ such that
$$\|Y\|_{{\s}^2}:=
  \left(\EX
    \left[
       \sup_{t\in\tT}|Y_t|^2
    \right]
  \right)^{1/2}<+\infty. $$
Moreover, let $\M^2(0,T;\rtn^{k\times d})$
denote the set of (equivalent classes of) $(\F_t)$-progressively measurable
${\rtn}^{k\times d}$-valued processes $(Z_t)_{t\in\tT}$ such that
$$\|Z\|_{{\M}^2}:=
  \left(\EX
    \left[
      \int_0^T |Z_t|^2\dif t
    \right]
  \right)^{1/2}<+\infty.$$
Obviously, $\s^2(0,T;\rtn^k)$ is a Banach space and $\M^2(0,T;\rtn^{k\times d})$
is a Hilbert space.

Finally, let $\BS$ be the set of all non-decreasing continuous functions
$\rho(\cdot):\rtn^+\mapsto\rtn^+$ with $\rho(0)=0$ and $\rho(x)>0$ for all $x>0$,
where and hereafter $\rtn^+:=[0,+\infty)$.

As mentioned above, we will deal only with the multidimensional BSDE which is an
equation of type \eqref{eq:BSDEs}, where the terminal condition $\xi$ is
$\F_T$-measurable, the terminal time $T$ satisfies $0\leq T\leq +\infty$, and the
generator $g$ is $(\F_t)$-progressively measurable for each $(y,z)$. In this
paper, we use the following definition.

\begin{dfn}
  A pair of processes $(y_t,z_t)_{t\in\tT}$ taking values in $\rtn^k\times\rtn^{k\times d}$
  is called a solution of BSDE \eqref{eq:BSDEs}, if $(y_t,z_t)_{t\in\tT}$ belongs
  to the space $\s^2(0,T;\rtn^k)\times\M^2(0,T;\rtn^{k\times d})$ and $\prs$,
  BSDE \eqref{eq:BSDEs} holds true for each $t\in\tT$.
\end{dfn}

The following Lemma \ref{lem:LepeltierMartinLemma} comes from \citet*{LepeltierSanMartin1997SPL},
which will be used later.

\begin{lem}\label{lem:LepeltierMartinLemma}
  Let $p\in\N$, $f(\cdot):\rtn^p\mapsto\rtn$ be a continuous and linear-growth function,
  i.e., there exists a positive constant $K$ such that
  $|f(x)|\leq K(1+|x|)$ for all
  $x\in\rtn^p$. Then $f_n(x)=\inf_{y\in\rtn^p}\{f(y)+n|x-y|\}$, $x\in\rtn^p$,
  is well defined for $n\geq K$ and satisfies
  \begin{enumerate}
    \renewcommand{\theenumi}{(\roman{enumi})}
    \renewcommand{\labelenumi}{\theenumi}
    \item Linear growth: for each $x\in\rtn^p$, $|f_n(x)|\leq K(1+|x|)$;
    \item Monotonicity in $n$: for each $x\in\rtn^p$, $f_n(x)$ increases in $n$;
    \item Lipschitz continuous: for each $x_1$, $x_2\in\rtn^p$, we have
           $|f_n(x_1)-f_n(x_2)|\leq n|x_1-x_2|$;
    \item Strong convergence: if $x_n\to x$, then $f_n(x_n)\to f(x)$ as $n\to+\infty$.
  \end{enumerate}
\end{lem}

In the following, we will establish some propositions on DBDEs
with general time intervals, which will play important roles in the proof of our
main result. It is very likely that these results have already appeared somewhere,
but we have not seen it, so we provide their proofs in Appendix for the
convenience of readers.

\begin{pro}\label{pro:DBDEResult}
  Let $0\leq T\leq+\infty$ and $f(t,y):\tT\times\rtn\mapsto\rtn$ satisfy the
  following two assumptions:
  \begin{enumerate}
    \renewcommand{\theenumi}{(B\arabic{enumi})}
    \renewcommand{\labelenumi}{\theenumi}
    \item \label{B:fNonuniformlyLipschitzContinuous}
          there exists a function $u(\cdot):\rtn^+\mapsto\rtn^+$ with
          $\intT{0}u(t)\dif t<+\infty$ such that for each $y_1$, $y_2\in\rtn$
          and $t\in\tT$,
          \begin{equation*}
            |f(t,y_1)-f(t,y_2)|\leq u(t)|y_1-y_2|;
          \end{equation*}
    \item \label{B:fOnlyIntegrable}
          $\intT{0}|f(t,0)|\dif t<+\infty$.
  \end{enumerate}
  Then for each $\delta\in\rtn$ the following DBDE
  \begin{equation}\label{eq:DBDEs}
    y_t=\delta+\intT{t}f(s,y_s)\dif s,\quad t\in\tT,
  \end{equation}
  has a unique continuous solution $(y_t)_{t\in\tT}$ such that
  $\sup_{t\in\tT}|y_t|<+\infty$.
\end{pro}

\begin{pro}\label{pro:DBDESolutionsDefinedRecursively}
  Assume $0\leq T\leq+\infty$, $f$ satisfies \ref{B:fNonuniformlyLipschitzContinuous}
  and \ref{B:fOnlyIntegrable}, $(y_t)_{t\in\tT}$ is the unique continuous solution of
  DBDE \eqref{eq:DBDEs} such that $\sup_{t\in\tT}|y_t|<+\infty$, $C$ is an
  arbitrary constant and $y^n_t$ is defined recursively as follows, for each $n\in\N$
  and $\delta\in\rtn$,
  \begin{equation*}
    y^1_t=C; \quad y^{n+1}_t=\delta+\intT{t}f(s,y^n_s)\dif s,
    \quad t\in\tT.
  \end{equation*}
  Then $y^n_t\to y_t$ as $n\to+\infty$ for each $t\in\tT$.
\end{pro}

\begin{pro}\label{pro:DBDESolutionComparison}
  Let $0\leq T\leq+\infty$, $f$ and $f'$ satisfy \ref{B:fNonuniformlyLipschitzContinuous}
  -- \ref{B:fOnlyIntegrable}, $(y_t)_{t\in\tT}$ and $(y'_t)_{t\in\tT}$ with $\sup_{t\in\tT}\big(|y_t|+|y'_t|\big)<+\infty$
  satisfy respectively DBDE \eqref{eq:DBDEs} and the following DBDE, for some $\delta'\in\rtn$,
  \begin{equation*}
    y'_t=\delta'+\intT{t}f'(s,y'_s)\dif s,\quad t\in\tT.
  \end{equation*}
  Assume that $f(t,y'_t)\geq f'(t,y'_t)$ for each $t\in\tT$. Then, we have
  \begin{enumerate}
    \renewcommand{\theenumi}{(\roman{enumi})}
    \renewcommand{\labelenumi}{\theenumi}
    \item (Comparison theorem) if $\delta\geq\delta'$, then $y_t\geq y'_t$ for each $t\in\tT$;
    \item (Strict comparison theorem) if $\delta>\delta'$, then $y_t>y'_t$
          for each $t\in\tT$.
  \end{enumerate}
\end{pro}

\begin{pro}\label{pro:DBDEsWithTwoFurthermoreResults}
  Let $0\leq T\leq+\infty$, $u(\cdot)$ be defined in \ref{B:fNonuniformlyLipschitzContinuous} and
  $\varphi(\cdot):\rtn^+\mapsto\rtn^+$ be a continuous function such that
  $\varphi(x)\leq ax+b$ for all $x\in\rtn^+$, where $a$ and $b$ are two given nonnegative
  constants. Then for each $\delta\in\rtn^+$, the following DBDE
  \begin{equation}\label{eq:DBDEFurthermoreResults}
    y^\delta_t=\delta+\intT{t}u(s)\varphi(y^\delta_s)\dif s,\quad t\in\tT,
  \end{equation}
  has a solution $(y^\delta_t)_{t\in\tT}$ such that $\sup_{t\in\tT}|y^\delta_t|<+\infty$.
  In addition,
  \begin{enumerate}
    \renewcommand{\theenumi}{(\roman{enumi})}
    \renewcommand{\labelenumi}{\theenumi}
    \item if $\delta>0$ and $\varphi(x)>0$ for all $x>0$, then DBDE \eqref{eq:DBDEFurthermoreResults}
          has a unique solution;
    \item if $\delta=0$ and $\varphi\in\BS$ with
          $\int_{0^+}\varphi^{-1}(x)\dif x=+\infty$, then DBDE \eqref{eq:DBDEFurthermoreResults}
          has a unique solution $y_t\equiv 0$.
  \end{enumerate}
\end{pro}
\section{Main result}
\label{sec:MainResultOnBSDEs}

In this section, we will state the main result of this paper.
Let us first introduce the following assumptions with respect to the
generator $g$ of BSDE \eqref{eq:BSDEs}, where $0\leq T\leq+\infty$.

\begin{enumerate}
  \renewcommand{\theenumi}{(H\arabic{enumi})}
  \renewcommand{\labelenumi}{\theenumi}
  \item \label{H:gUniformlyContinuousInY}
        $g$ is uniformly continuous in $y$ non-uniformly with respect to $t$,
        i.e., there exists a deterministic function $u(\cdot):\tT\mapsto\rtn^+$
        with $\intT{0}u(t)\dif t<+\infty$ and a linear-growth function
        $\rho(\cdot)\in\BS$ such that $\pts$,
        for each $y_1$, $y_2\in\rtn^k$ and $z\in\rtn^{k\times d}$,
        \begin{equation*}
          |g(\omega,t,y_1,z)-g(\omega,t,y_2,z)|\leq u(t)\rho(|y_1-y_2|);
        \end{equation*}
        Furthermore, we also assume that $\int_{0^+}\rho^{-1}(u)\dif u=+\infty$;
  \item \label{H:gUniformlyContinuousInZ}
        $g$ is uniformly continuous in $z$ non-uniformly with respect to $t$,
        i.e., there exists a deterministic function $v(\cdot):\tT\mapsto\rtn^+$ with
        $\intT{0}\big(v(t)+v^2(t)\big)\dif t<+\infty$ and a linear-growth function
        $\phi(\cdot)\in\BS$
        such that $\pts$, for each $y\in\rtn^k$ and $z_1$, $z_2\in\rtn^{k\times d}$,
        \begin{equation*}
          |g(\omega,t,y,z_1)-g(\omega,t,y,z_2)|\leq v(t)\phi(|z_1-z_2|);
        \end{equation*}
  \item \label{H:ithComponent}
        For any $i=1,\cdots,k$, $g_i(\omega,t,y,z)$, the $i$th component of $g$, depends
        only on the $i$th row of $z$;
  \item \label{H:gSquareIntegrable}
        $\EXlr{\left(\intT{0}|g(\omega,t,0,0)|\dif t\right)^2}<+\infty$.
\end{enumerate}

In the sequel, we denote the linear-growth constant for $\rho(\cdot)$ and
$\phi(\cdot)$ in \ref{H:gUniformlyContinuousInY} and \ref{H:gUniformlyContinuousInZ}
by $A>0$, i.e., $\rho(x)\leq A(1+x)$ and $\phi(x)\leq A(1+x)$ for all $x\in\rtn^+$.
In the remaining of this paper, we put an $i$ at upper left of
$y\in\rtn^k$, $z\in\rtn^{k\times d}$ to represent the $i$th component of $y$ and
the $i$th row of $z$, like ${}^iy$ and ${}^iz$.

The main result of this paper is the following Theorem \ref{thm:MainResult}, whose
proof will be given in next section.

\begin{thm}\label{thm:MainResult}
Assume that $0\leq T\leq+\infty$ and $g$ satisfies \ref{H:gUniformlyContinuousInY}
-- \ref{H:gSquareIntegrable}. Then for each $\xi\in L^2(\Omega,\F_T,\PR;\rtn^k)$,
BSDE \eqref{eq:BSDEs} has a unique solution.
\end{thm}

\begin{rmk}\label{rmk:GeneralizeResultInHamadene}
In the corresponding assumptions in \citet*{Hamadene2003Bernoulli} and
\citet*{FanJiangDavison2010CRASSI} the $u(t)$,
$v(t)$ appearing in \ref{H:gUniformlyContinuousInY} and \ref{H:gUniformlyContinuousInZ}
are bounded by a constant $c>0$ , and $T$ is a finite real number.
However, in our framework the $u(t)$, $v(t)$ may be unbounded. In addition, Theorem
\ref{thm:MainResult} also considers the case of $T=+\infty$. Consequently,
Theorem \ref{thm:MainResult} generalizes the corresponding results in
\citet*{Hamadene2003Bernoulli} and \citet*{FanJiangDavison2010CRASSI}.
\end{rmk}

\begin{ex}
  Let $0\leq T\leq +\infty$, and for each $i=1,\cdots,k$ and $(\omega,t,y,z)\in\Omega\times\tT\times\rtn^k\times\rtn^{k\times d}$,
  define the generator $g=(g_1,\cdots,g_k)$ by
  \begin{equation*}
    g_i(\omega,t,y,z)=f_1(t)\big(h(|y|)+1\big)+f_2(t)\sqrt{|{}^iz|}+|B_t(\omega)|,
  \end{equation*}
  where
  \begin{align*}
    f_1(t)&=\frac{1}{\sqrt{t}}\one{0<t<\delta}+\frac{1}{\sqrt{1+t^2}}\one{\delta\leq t\leq T},\\
    f_2(t)&=\frac{1}{\sqrt[4]{t}}\one{0<t<\delta}+\frac{1}{\sqrt{(1+t)^2}}\one{\delta\leq t\leq T},\\
    h(x)&=x\ln\frac{1}{x}\one{0\leq x\leq\delta}+[h'(\delta-)(x-\delta)+h(\delta)]\one{x>\delta},
  \end{align*}
  with $\delta$ small enough. Since $h(0)=0$ and $h$ is concave and increasing,
  we have $h(x_1+x_2)\leq h(x_1)+h(x_2)$ for all $x_1$, $x_2\in\rtn^+$, which
  implies that $|h(x_1)-h(x_2)|\leq h(|x_1-x_2|)$. Thus, note that
  $$\int_{0^+}\frac{1}{x\ln\frac{1}{x}}\dif x=+\infty.$$
  We know that the generator $g$ satisfies assumptions
  \ref{H:gUniformlyContinuousInY} -- \ref{H:gSquareIntegrable}
  with $u(t)=f_1(t)$, $v(t)=f_2(t)$. It then follows from Theorem \ref{thm:MainResult}
  that for each $\xi\in L^2(\Omega,\F_T,\PR;\rtn^k)$,
  BSDE \eqref{eq:BSDEs} has a unique solution $(y_t,z_t)_{t\in\tT}$.

  It should be mentioned that the above conclusion can not be obtained by the result of
  \citet*{Hamadene2003Bernoulli}, \citet*{FanJiangDavison2010CRASSI} and other existing results.
\end{ex}

\section{Proof of the main result}
\label{sec:ProofOfTheorem}
This section will give the proof of our main result --- Theorem \ref{thm:MainResult}.
Before starting the proof, let us first introduce the following Lemma
\ref{lem:ResultInChenWang}, which comes from Theorem 1.2 in \citet*{ChenWang2000JAMSA}.
The following assumption will be used in Lemma \ref{lem:ResultInChenWang}, where we suppose $0\leq T\leq+\infty$:

\begin{enumerate}
  \renewcommand{\theenumi}{(A\arabic{enumi})}
  \renewcommand{\labelenumi}{\theenumi}
  \item \label{A:gNonuniformlyLipschitzContinuousInYAndZ}
        There exist two deterministic functions $u(\cdot),v(\cdot): \tT\mapsto\rtn^+$
        with $\intT{0}\big(u(t)+v^2(t)\big)\dif t<+\infty$ such that $\pts$, for each
        $y_1$, $y_2\in\rtn^k$ and $z_1$, $z_2\in\rtn^{k\times d}$,
        \begin{equation*}
          |g(t,y_1,z_1)-g(t,y_2,z_2)|\leq u(t)|y_1-y_2|+v(t)|z_1-z_2|.
        \end{equation*}
\end{enumerate}

\begin{lem}[Theorem 1.2 in \citet*{ChenWang2000JAMSA}]\label{lem:ResultInChenWang}
  Assume that $0\leq T\leq+\infty$ and $g$ satisfies \ref{A:gNonuniformlyLipschitzContinuousInYAndZ}
  and \ref{H:gSquareIntegrable}. Then for each $\xi\in L^2(\Omega,\F_T,\PR;\rtn^k)$,
  BSDE \eqref{eq:BSDEs} has a unique solution $(y_t,z_t)_{t\in\tT}$.
\end{lem}

\subsection{Proof of the uniqueness part of Theorem \ref{thm:MainResult}}

  The idea of the proof of this part is partly motivated by \citet*{FanJiangDavison2010CRASSI}.
  Let $(y^1_t,z^1_t)_{t\in\tT}$ and $(y^2_t,z^2_t)_{t\in\tT}$ be two solutions of BSDE
  \eqref{eq:BSDEs}. Then
  we have the following Lemma \ref{lem:YiYiiUniformlyBoundedAndConditionalExpectation},
  whose proof is provided at the end of this subsection.
  \begin{lem}\label{lem:YiYiiUniformlyBoundedAndConditionalExpectation}
    The process $(y^1_t-y^2_t)_{t\in\tT}$ is uniformly bounded, i.e., there
    exists a positive constant $C_1>0$ such that
    \begin{equation}\label{eq:YiAndYiiUniformlyBoundedCi}
    \pts,\quad |y^1_t-y^2_t|\leq C_1.
    \end{equation}
    Moreover, for each $n\in\N$, $i=1,2,\cdots,k$ and $0\leq r\leq t\leq T$, we have
    \begin{equation}\label{eq:ExpectationMAndIYiAndYiiIthComponentLeqInequality}
      \EX^{n,i}
      \left[
        |{}^iy^1_t-{}^iy^2_t|\big|\F_r
      \right]\leq a_n+
      \intT{t}\EX^{n,i}
      \left[
        u(s)\rho(|y^1_s-y^2_s|)\big|\F_r
      \right]\dif s,
    \end{equation}
    where
    $$a_n=\phi\left(\frac{2A}{n+2A}\right)\intT{0}v(s)\dif s,$$
    and $\EX^{n,i}[X|\F_t]$ represents the conditional expectation of random
    variable $X$ with respect to $\F_t$ under a probability measure
    $\PR^{n,i}$ on $(\Omega,\F)$, which depends on $n$ and $i$, and which is absolutely
    continuous with respect to $\PR$.
  \end{lem}

  In the sequel, let $\overline{\rho}(y)=\rho(|y|)$ for each $y\in\rtn$, and
  for each $n\in\N$, define $\rho_{n}(\cdot):\rtn\mapsto\rtn^+$ by
  $$\rho_{n}(x)=\sup_{y\in\rtn}\{\overline\rho(y)-n|x-y|\}.$$
  It follows from Lemma \ref{lem:LepeltierMartinLemma} that $\rho_{n}$ is well
  defined for $n\geq A$, Lipschitz continuous, non-increasing in $n$ and converges
  to $\overline\rho$. Then, for each $n\geq A$, by Proposition \ref{pro:DBDEResult} we can
  let $f^{n}_t$ be the unique solution of the following DBDE
  \begin{equation}\label{eq:UniquenessPartDBDEWithRhoN}
    f^{n}_t=a_n+\intT{t}[u(s)\rho_n(k\cdot f^{n}_s)]\dif s,\quad t\in\tT.
  \end{equation}

  Noticing that $\rho_{n}$ and $a_n$ are both decreasing in $n$, we have
  $0\leq f^{n+1}_t\leq f^{n}_t$ for each $n\geq A$ by Proposition
  \ref{pro:DBDESolutionComparison}, which implies that the sequence
  $\{f^{n}_t\}^{+\infty}_{n=1}$
  converges point wisely to a function $f_t$. Thus, by sending $n\to+\infty$ in
  \eqref{eq:UniquenessPartDBDEWithRhoN}, it follows from Lemma \ref{lem:LepeltierMartinLemma}
  and the Lebesgue dominated convergence theorem that
  $$f_t=\intT{t}[u(s)\overline\rho(k\cdot f_s)]\dif s
       =\intT{t}[u(s)\rho(k\cdot f_s)]\dif s,\quad t\in\tT.$$
  Recalling that $\rho(\cdot)\in\BS$ and $\int_{0^+}\rho^{-1}(u)\dif u=+\infty$,
  Proposition \ref{pro:DBDEsWithTwoFurthermoreResults} yields that $f_t\equiv 0$.

  Now, for each $n\geq A$, $j\in\N$ and $t\in\tT$, let $f^{n,j}_t$ be the function defined
  recursively as follows:
  \begin{equation}\label{eq:FNAndJRecursivelyDefined}
    f^{n,1}_t=C_1;\qquad
    f^{n,j+1}_t=a_n
      +\intT{t}\big(u(s)\rho_{n}(k\cdot f^{n,j}_s)\big)\dif s,
  \end{equation}
  where $C_1$ is defined in \eqref{eq:YiAndYiiUniformlyBoundedCi}.
  Noticing that $\rho_{n}$ is Lipschitz continuous, by Proposition
  \ref{pro:DBDESolutionsDefinedRecursively} we know that $f^{n,j}_t$ converges point wisely
  to $f^n_t$ as $j\to+\infty$ for each $t\in\tT$ and $n\geq A$.

  On the other hand, it is easy to check by induction that for each $n\geq A$,
  $j\geq 1$ and $i=1,\cdots,k$,
  \begin{equation}\label{eq:YiAndYiiIthComponentControlByFNAndJt}
    |{}^iy^1_t-{}^iy^2_t|\leq f^{n,j}_t\leq f_0^{n,j},\quad t\in\tT.
  \end{equation}
  Indeed, \eqref{eq:YiAndYiiIthComponentControlByFNAndJt} holds true for $j=1$ due
  to \eqref{eq:YiAndYiiUniformlyBoundedCi}. Suppose
  \eqref{eq:YiAndYiiIthComponentControlByFNAndJt} holds true for $j\geq 1$. Then,
  for each $t\in\tT$,
  \begin{equation*}
    u(t)\rho(|y^1_t-y^2_t|)
    \leq u(t)\rho(k\cdot f^{n,j}_t)
    \leq u(t)\rho_n(k\cdot f^{n,j}_t).
  \end{equation*}
  In view of \eqref{eq:ExpectationMAndIYiAndYiiIthComponentLeqInequality} with
  $r=t$ as well as \eqref{eq:FNAndJRecursivelyDefined}, we can deduce that for each
  $n\geq A$ and $i=1,2,\cdots,k$,
  \begin{equation*}
    |{}^iy^1_t-{}^iy^2_t|\leq f^{n,j+1}_t\leq f^{n,j+1}_0,\quad t\in\tT,
  \end{equation*}
  which is the desired result.

  Finally, by sending first $j\to+\infty$ and then $n\to+\infty$
  in \eqref{eq:YiAndYiiIthComponentControlByFNAndJt}, we obtain that
  $\sup_{t\in\tT}|{}^iy^1_t-{}^iy^2_t|=0$ for each $i=1,2,\cdots,k$. That is,
  the solution of BSDE \eqref{eq:BSDEs} is unique. The proof of the uniqueness
  part is then completed.\hfill$\square$

\begin{proof}[{\bf Proof of Lemma \ref{lem:YiYiiUniformlyBoundedAndConditionalExpectation}}]
  Using \Ito{} formula to $|y^1_t-y^2_t|^2$ we arrive, for each $t\in\tT$, at
  \begin{align}
    |y^1_t-y^2_t|^2+\intT{t}|z^1_s-z^2_s|^2\dif s
    =&\ 2\intT{t}\langle y^1_s-y^2_s,g(s,y^1_s,z^1_s)-g(s,y^2_s,z^2_s)\rangle\dif s\nonumber\\
     &\ -2\intT{t}\langle y^1_s-y^2_s,(z^1_s-z^2_s)\dif B_s\rangle.\label{eq:AfterItoFormulaYiAndYiiUniformlyBounded}
  \end{align}
  The inner product term including $g$ can be enlarged by
  \ref{H:gUniformlyContinuousInY} -- \ref{H:gUniformlyContinuousInZ} and the
  basic inequality $2ab\leq 2a^2+b^2/2$ as follows:
  \begin{align*}
    2\langle y^1_s-y^2_s,g(s,y^1_s,z^1_s)-g(s,y^2_s,z^2_s)\rangle
    & \leq 2|y^1_s-y^2_s||g(s,y^1_s,z^1_s)-g(s,y^2_s,z^1_s)+g(s,y^2_s,z^1_s)
                         -g(s,y^2_s,z^2_s)|\\
    & \leq 2|y^1_s-y^2_s|\big(Au(s)|y^1_s-y^2_s|+Av(s)|z^1_s-z^2_s|+Au(s)+Av(s)\big)\\
    & \leq B(s)|y^1_s-y^2_s|^2+\frac{1}{2}|z^1_s-z^2_s|^2+Au(s)+Av(s),
  \end{align*}
  where $B(s)=2Au(s)+2A^2v^2(s)+Au(s)+Av(s)$.
  Putting the previous inequality into \eqref{eq:AfterItoFormulaYiAndYiiUniformlyBounded}
  we can obtain that for each $t\in\tT$,
  \begin{equation*}
    |y^1_t-y^2_t|^2\leq\intT{t}B(s)|y^1_s-y^2_s|^2\dif s
    -2\intT{t}\langle y^1_s-y^2_s,(z^1_s-z^2_s)\dif B_s\rangle+C,
  \end{equation*}
  where $C=\intT{0}A(u(s)+v(s))\dif s$. Note that both $(y^1_t,z^1_t)_{t\in\tT}$ and
  $(y^2_t,z^2_t)_{t\in\tT}$ belong to the process space $\s^2(0,T;\rtn^k)\times\M^2(0,T;\rtn^{k\times d})$.
  By the Burkholder-Davis-Gundy (BDG for short in the remaining) inequality and
  H\"older's inequality we have that there exists a positive constant $K'>0$
  such that
  \begin{align*}
    \EXlr{\sup_{t\in\tT}
    \left|
      \intT[t]{0}\langle y^1_s-y^2_s,(z^1_s-z^2_s)\dif B_s
    \right|}
    & \leq K'\EXlr{\sqrt{\intT{0}|y^1_s-y^2_s|^2|z^1_s-z^2_s|^2\dif s}}\\
    & \leq K'\sqrt{\EXlr{\sup_{t\in\tT}|y^1_t-y^2_t|^2}}
             \sqrt{\EXlr{\intT{0}|z^1_s-z^2_s|^2\dif s}}<+\infty,
  \end{align*}
  which implies that $(\intT[t]{0}\langle y^1_s-y^2_s,(z^1_s-z^2_s)\dif B_s\rangle)_{t\in\tT}$
  is an $(\F_t,\PR)$-martingale.
  Then we have for each $0\leq r\leq t\leq T$,
  \begin{equation*}
    \EXlr{|y^1_t-y^2_t|^2\big|\F_r}\leq
    \intT{t}B(s)\EXlr{|y^1_s-y^2_s|^2\big|\F_r}\dif s+C.
  \end{equation*}
  By Lemma 4 in \citet*{FanJiangTian2011SPA} we have
  $$\EXlr{|y^1_t-y^2_t|^2\big|\F_r}\leq C\me^{\intT{0}B(s)\dif s}:=(C_1)^2,$$
  which yields \eqref{eq:YiAndYiiUniformlyBoundedCi} after taking $r=t$.

  In the sequel, by \ref{H:ithComponent} we have for each $t\in\tT$,
  \begin{equation*}%\label{eq:BSDEYiYiiIthComponent}
    {}^iy^1_t-{}^iy^2_t=\intT{t}\big(g_i(s,y^1_s,{}^iz^1_s)-g_i(s,y^2_s,{}^iz^2_s)\big)\dif s
    -\intT{t}({}^iz^1_s-{}^iz^2_s)\dif B_s,
  \end{equation*}
  Then, \ref{H:ithComponent} and Tanaka's formula lead to that, for each $t\in\tT$,
  \begin{equation}\label{eq:YiYiiIthComponentAfterTanakaFormula}
    |{}^iy^1_t-{}^iy^2_t|
    \leq\intT{t}\sgn({}^iy^1_s-{}^iy^2_s)\big(g_i(s,y^1_s,{}^iz^1_s)-g_i(s,y^2_s,{}^iz^2_s)\big)\dif s
    -\intT{t}\sgn({}^iy^1_s-{}^iy^2_s)({}^iz^1_s-{}^iz^2_s)\dif B_s.
  \end{equation}
  Furthermore, it follows from \ref{H:gUniformlyContinuousInY} and \ref{H:gUniformlyContinuousInZ}
  that
  \begin{equation}\label{eq:GIthComponentControl}
    |g_i(s,y^1_s,{}^iz^1_s)-g_i(s,y^2_s,{}^iz^2_s)|
    \leq u(s)\rho(|y^1_s-y^2_s|)+v(s)\phi(|{}^iz^1_s-{}^iz^2_s|).
  \end{equation}
  Recalling that $\phi(\cdot)$ is a non-decreasing function from $\rtn^+$
  to itself with at most linear-growth. From \citet*{FanJiangDavison2010CRASSI}
  we know that for each $n\in\N$ and $x\in\rtn^+$,
  \begin{equation}\label{eq:PhiControlAndDivide}
    \phi(x)\leq(n+2A)x+\phi\left(\frac{2A}{n+2A}\right).
  \end{equation}
  Thus, combining \eqref{eq:YiYiiIthComponentAfterTanakaFormula}
  -- \eqref{eq:PhiControlAndDivide} we get that for each $n\in\N$,
  \begin{align*}
    |{}^iy^1_t-{}^iy^2_t|
    \leq&\ \phi\left(\frac{2A}{n+2A}\right)\intT{0}v(s)\dif s
       +\intT{t}\big(u(s)\rho(|y^1_s-y^2_s|)+(n+2A)v(s)|{}^iz^1_s-{}^iz^2_s|\big)\dif s\\
    &\ -\intT{t}\sgn({}^iy^1_s-{}^iy^2_s)({}^iz^1_s-{}^iz^2_s)\dif B_s,\quad t\in\tT.
  \end{align*}
  Now for each $t\in\tT$, let
  \begin{equation*}
    e^{n,i}_t:=(n+2A)
    \frac{\sgn({}^iy^1_t-{}^iy^2_t)({}^iz^1_t-{}^iz^2_t)^*}{|{}^iz^1_t-{}^iz^2_t|}
    \one{|{}^iz^1_t-{}^iz^2_t|\neq 0}.
  \end{equation*}
  Then, $(e^{n,i}_t)_{t\in\tT}$ is a $\rtn^d$-valued, bounded and
  $(\F_t)$-adapted process. It follows from Girsanov's theorem that
  $B^{n,i}_t=B_t-\intT[t]{0}e^{n,i}_sv(s)\dif s$,
  $t\in\tT$, is a $d$-dimensional Brownian motion under the probability
  $\PR^{n,i}$ on $(\Omega,\F)$ defined by
  \begin{equation*}
    \frac{\dif\PR^{n,i}}{\dif\PR}
    =\exp\left\{\intT{0}v(s)(e^{n,i}_s)^*\dif B_s
      -\frac{1}{2}\intT{0}v^2(s)|e^{n,i}_s|^2\dif s\right\}.
  \end{equation*}
  Thus, for each $n\in\N$ and $t\in\tT$,
  \begin{equation}\label{eq:YiAndYiiIComponentBeforeBDGInequality}
    |{}^iy^1_t-{}^iy^2_t|
    \leq\phi\left(\frac{2A}{n+2A}\right)
            \intT{0}v(s)\dif s
    +\intT{t}u(s)\rho(|y^1_s-y^2_s|)\dif s
    -\intT{t}\sgn({}^iy^1_s-{}^iy^2_s)({}^iz^1_s-{}^iz^2_s)\dif B^{n,i}_s.
  \end{equation}
  Moreover, the process
  $\big(\intT[t]{0}\sgn({}^iy^1_s-{}^iy^2_s)({}^iz^1_s-{}^iz^2_s)
   \dif B^{n,i}_s\big)_{t\in\tT}$ is an $(\F_t,\PR^{n,i})$-martingale.
  In fact, let $\EX^{n,i}[X]$ represent the expectation of the
  random variable $X$ under $\PR^{n,i}$. By the BDG
  inequality and H\"older's inequality we know that there exists
  a positive constant $K''>0$ such that for each $n\in\N$,
   \begin{align*}
     & \EX^{n,i}
       \left[
         \sup_{t\in\tT}
         \left|
           \intT[t]{0}\sgn({}^iy^1_s-{}^iy^2_s)({}^iz^1_s-{}^iz^2_s)
           \dif B^{n,i}_s
         \right|
       \right]\\
     & \leq K''\EX^{n,i}
       \left[
         \sqrt{\intT{0}|{}^iz^1_s-{}^iz^2_s|^2\dif s}
       \right]
       \leq K''\sqrt{\EXlr{\left(\frac{\dif\PR^{n,i}}{\dif\PR}\right)^2}}
            \sqrt{\EXlr{\intT{0}|{}^iz^1_s-{}^iz^2_s|^2\dif s}}<+\infty.
   \end{align*}
   Thus, for each $n\in\N$ and $0\leq r\leq t\leq T$, by taking the condition
   expectation with respect to $\F_r$ under $\PR^{n,i}$ in both sides of
   \eqref{eq:YiAndYiiIComponentBeforeBDGInequality}, we can get the desired
   result \eqref{eq:ExpectationMAndIYiAndYiiIthComponentLeqInequality}.
   The proof of Lemma \ref{lem:YiYiiUniformlyBoundedAndConditionalExpectation}
   is complete.
\end{proof}

\subsection{Proof of the existence part of Theorem \ref{thm:MainResult}}

  The idea of the proof of this part is enlightened by \citet*{Hamadene2003Bernoulli}.
  But some different arguments are used, and then the proof procedure is simplified at certain degree.

  Assume that the generator $g$ satisfies \ref{H:gUniformlyContinuousInY} -- \ref{H:gSquareIntegrable}
  and $\xi\in L^2(\Omega,\F_T,\PR;\rtn^k)$. Without loss of generality, we assume
  that $u(t)$ and $v(t)$ in \ref{H:gUniformlyContinuousInY} and \ref{H:gUniformlyContinuousInZ}
  are both strictly positive functions. Otherwise, we can use $u(t)+\me^{-t}$ and
  $v(t)+\me^{-t}$ instead of them respectively.

  We have the following lemma
  whose proof is placed at the end of this subsection.

\begin{lem}\label{lem:GnSequenceConvolution}
  Let $g$ satisfy \ref{H:gUniformlyContinuousInY} -- \ref{H:ithComponent}, and
  assume that $u(t)>0$ and $v(t)>0$ for each $t\in\tT$.
  Then there exists a generator sequence $\{g^n\}^{+\infty}_{n=1}$ such that
  \begin{enumerate}
    \renewcommand{\theenumi}{(\roman{enumi})}
    \renewcommand{\labelenumi}{\theenumi}
    \item For each $n\in\N$, $g^n(t,y,z)$ is a mapping from
          $\Omega\times\tT\times\rtn^k\times\rtn^{k\times d}$ into $\rtn^k$ and is
          $(\F_t)$-progressively measurable. Moreover, we have $\pts$, for each
          $y\in\rtn^k$ and $z\in\rtn^{k\times d}$,
          \begin{equation*}
            |g^n(t,y,z)|\leq |g(t,0,0)|+kAu(t)(1+|y|)+kAv(t)(1+|z|);
          \end{equation*}
    \item For each $n\in\N$, $g^n(t,y,z)$ satisfies \ref{H:ithComponent}, and
          $\pts$, for each $y_1$, $y_2\in\rtn^k$ and $z_1$, $z_2\in\rtn^{k\times d}$,
          we have
          \begin{align*}
            |g^n(t,y_1,z_1)-g^n(t,y_2,z_2)|
            & \leq ku(t)\rho(|y_1-y_2|)+kv(t)\phi(|z_1-z_2|),\\
            |g^n(t,y_1,z_1)-g^n(t,y_2,z_2)|
            & \leq k(n+A)\big(u(t)|y_1-y_2|+v(t)|z_1-z_2|\big);
          \end{align*}
    \item For each $n\in\N$, there exists a non-increasing deterministic functions sequence
          $b_n(\cdot):\tT\mapsto\rtn^+$ with $\intT{0}b_n(t)\dif t\to 0$
          as $n\to +\infty$ such that $\pts$, for each $y\in\rtn^k$ and
          $z\in\rtn^{k\times d}$,
          \begin{equation*}
            |g^n(t,y,z)-g(t,y,z)| \leq kb_n(t).
          \end{equation*}
  \end{enumerate}
\end{lem}

  It follows from (i) -- (ii) of Lemma \ref{lem:GnSequenceConvolution} and \ref{H:gSquareIntegrable}
  that for each $n\in\N$, $g^n$ satisfies \ref{A:gNonuniformlyLipschitzContinuousInYAndZ}
  and \ref{H:gSquareIntegrable}. Then it follows from Lemma \ref{lem:ResultInChenWang}
  that for each $n\in\N$ and $\xi\in L^2(\Omega,\F_T,\PR;\rtn^k)$, the following BSDE
  \begin{equation}\label{eq:BSDEsYnInExistenceProof}
    y^n_t=\xi+\intT{t}g^n(s,y^n_s,z^n_s)\dif s-\intT{t}z^n_s\dif B_s,\quad t\in\tT,
  \end{equation}
  has a unique solution $(y^n_t,z^n_t)_{t\in\tT}$.
  The following proof will be split into three steps.

  \textbf{Step 1.} In this step we show that $\{(y^n_t)_{t\in\tT}\}^{+\infty}_{n=1}$ is a
  Cauchy sequence in $\s^2(0,T;\rtn^k)$.

  For each $n$, $m\in\N$, let $(y^n_t,z^n_t)_{t\in\tT}$ and $(y^m_t,z^m_t)_{t\in\tT}$
  be, respectively, solutions of BSDE $(\xi,T,g^n)$ and BSDE $(\xi,T,g^m)$.
  Using \Ito{} formula to $|y^n_t-y^m_t|^2$ we arrive, for each $t\in\tT$, at
  \begin{align}
    |y^n_t-y^m_t|^2+\intT{t}|z^n_s-z^m_s|^2\dif s
    =&\ 2\intT{t}\langle y^n_s-y^m_s,g^n(s,y^n_s,z^n_s)-g^m(s,y^m_s,z^m_s)\rangle\dif s\nonumber\\
     &\ -2\intT{t}\langle y^n_s-y^m_s,(z^n_s-z^m_s)\dif B_s\rangle.\label{eq:YnAndYmAfterItoFormulaExistencePart}
  \end{align}
  It follows from (ii) -- (iii) in Lemma \ref{lem:GnSequenceConvolution}
  and the basic inequality $2ab\leq 2a^2+b^2/2$ that,
  with adding and subtracting the term $g^n(s,y^m_s,z^m_s)$,
  \begin{align*}
    & 2\langle y^n_s-y^m_s,g^n(s,y^n_s,z^n_s)-g^m(s,y^m_s,z^m_s)\rangle\\
    & \leq 2|y^n_s-y^m_s|\big(kAu(s)|y^n_s-y^m_s|+kAv(s)|z^n_s-z^m_s|+kAu(s)+kAv(s)+\tau_{n,m}(s)\big)\\
    & \leq \big(D(s)+\tau_{n,m}(s)\big)|y^n_s-y^m_s|^2+\frac{1}{2}|z^n_s-z^m_s|^2+kAu(s)+kAv(s)
      +\tau_{n,m}(s),
  \end{align*}
  where
  \begin{equation*}
    \tau_{n,m}(s)=k\big(b_n(s)+b_m(s)\big),\quad
    D(s)=2kAu(s)+2k^2A^2v^2(s)+kAu(s)+kAv(s).
  \end{equation*}
  Putting the previous inequality into \eqref{eq:YnAndYmAfterItoFormulaExistencePart}
  and taking the conditional expectation with respect to $\F_r$
  yield that, for each $0\leq r\leq t\leq T$ and $n$, $m\in\N$.
  \begin{equation*}
    \EXlr{|y^n_t-y^m_t|^2\big|\F_r}\leq
    \intT{t}\big(D(s)+\tau_{n,m}(s)\big)\EXlr{|y^n_s-y^m_s|^2\big|\F_r}\dif s+
    C_{n,m},
  \end{equation*}
  where $C_{n,m}=\intT{0}\big(kAu(s)+kAv(s)+\tau_{n,m}(s)\big)\dif s$.
  It follows from Lemma 4 in \citet*{FanJiangTian2011SPA} that
  $$\EXlr{|y^n_t-y^m_t|^2\big|\F_r}\leq C_{n,m}
    \me^{\intT{0}(D(s)+\tau_{n,m}(s))\dif s}
    \leq C_{1,1}\me^{\intT{0}(D(s)+\tau_{1,1}(s))\dif s}:=(C_2)^2.$$
%  Noticing that $\intT{0}\tau_{n,m}(s)\dif s\to 0$ as $n$, $m\to+\infty$.
  After taking $r=t$ in the previous inequality, we have that for each $n$, $m\in\N$,
  %\begin{equation*}%\label{eq:YnAndYmUniformlyBoundedCii}
    $\pts$, $|y^n_t-y^m_t|\!\leq\! C_2$.
  %\end{equation*}

  Furthermore, it follows from (ii) in Lemma \ref{lem:GnSequenceConvolution},
  \ref{H:ithComponent} and Tanaka's formula that for
  each $t\in\tT$,
  \begin{align}
    |{}^iy^n_t-{}^iy^m_t|
    \leq&\ \intT{t}\sgn({}^iy^n_s-{}^iy^m_s)\big(g^n_i(s,y^n_s,{}^iz^n_s)-g^m_i(s,y^m_s,{}^iz^m_s)\big)\dif s\nonumber\\
        &\ -\intT{t}\sgn({}^iy^n_s-{}^iy^m_s)({}^iz^n_s-{}^iz^m_s)\dif B_s.
        \label{eq:YnYmIthComponentAfterItoFormula}
  \end{align}
  It follows from (ii) -- (iii) in Lemma \ref{lem:GnSequenceConvolution} that, by
  adding and subtracting the term $g^n_i(s,y^m_s,{}^iz^m_s)$,
  \begin{equation}\label{eq:GnGmIthComponentControl}
    |g^n_i(s,y^n_s,{}^iz^n_s)-g_i^m(s,y^m_s,{}^iz^m_s)|
    \leq ku(s)\rho(|y^n_s-y^m_s|)+kv(s)\phi(|{}^iz^n_s-{}^iz^m_s|)+\tau_{n,m}(s).
  \end{equation}
  Combining \eqref{eq:YnYmIthComponentAfterItoFormula} -- \eqref{eq:GnGmIthComponentControl}
  with \eqref{eq:PhiControlAndDivide}
  we get that for each $n$, $m$, $q\in\N$ and $t\in\tT$,
  \begin{align}
    |{}^iy^n_t-{}^iy^m_t|
    \leq&\ C_{n,m,q}+k\intT{t}\big(u(s)\rho(|y^n_s-y^m_s|)+(q+2A)v(s)|{}^iz^n_s-{}^iz^m_s|\big)\dif s\nonumber\\
    &\ -\intT{t}\sgn({}^iy^n_s-{}^iy^m_s)({}^iz^n_s-{}^iz^m_s)\dif B_s,
     \label{eq:YnAndYmIthComponentBeforeGirsanovTheorem}
  \end{align}
  where
  $$C_{n,m,q}:=k\phi
    \left(
      \frac{2A}{q+2A}
    \right)\intT{0}v(s)\dif s+\intT{0}\tau_{n,m}(s)\dif s.$$
  In the sequel, by virtue of Girsanov's theorem, in the same way as in the proof of Lemma
  \ref{lem:YiYiiUniformlyBoundedAndConditionalExpectation} we can deduce from
  \eqref{eq:YnAndYmIthComponentBeforeGirsanovTheorem} that for each $n$, $m$,
  $q\in\N$, $i=1,\cdots,k$, and $0\leq r\leq t\leq T$,
  \begin{equation*}%\label{eq:ExpectationMAndIYnAndYmIthComponentLeqInequality}
    \EX^{n,m,q,i}
      \left[
        |{}^iy^n_t-{}^iy^m_t|\big|\F_r
      \right]
      \leq C_{n,m,q}
        +k\intT{t}\EX^{n,m,q,i}
         \left[
           u(s)\rho(|y^n_s-y^m_s|)\big|\F_r
         \right]\dif s,
  \end{equation*}
  where $\EX^{n,m,q,i}[X|\F_t]$ represents the conditional expectation of random
  variable $X$ with respect to $\F_t$ under a probability measure
  $\PR^{n,m,q,i}$ on $(\Omega,\F)$, which depends on $n$, $m$, $q$ and $i$, and which
  is absolutely continuous with respect to $\PR$.

  Finally, note that $C_{n,m,q}$ tends non-increasingly to $0$ as $n$, $m$, $q\to+\infty$.
  The same argument as in the proof of the uniqueness part of Theorem
  \ref{thm:MainResult} yields that for each $i=1,\cdots,k$,
  \begin{equation*}
    \lim_{n,m\to+\infty} \EXlr{\sup_{t\in\tT}|{}^iy^n_t-{}^iy^m_t|^2}=0,
  \end{equation*}
  which means that $\{(y^n_t)_{t\in\tT}\}^{+\infty}_{n=1}$
  is a Cauchy sequence in $\s^2(0,T;\rtn^k)$. We denote the limit by
  $(y_t)_{t\in\tT}$.

\textbf{Step 2.} In this step we show that $\{(z^n_t)_{t\in\tT}\}^{+\infty}_{n=1}$ is
a Cauchy sequence in $\M^2(0,T;\rtn^{k\times d})$.

Using \Ito{} formula for $|y^n_t|^2$ defined in BSDE \eqref{eq:BSDEsYnInExistenceProof},
we can obtain that
\begin{equation*}
  |y^n_t|^2+\intT{t}|z^n_s|^2\dif s=|\xi|^2
  +2\intT{t}\langle y^n_s,g^n(s,y^n_s,z^n_s)\rangle\dif s-2\intT{t}\langle y^n_s,z^n_s\dif B_s\rangle.
\end{equation*}
Let $G_n(\omega):=\sup_{t\in\tT}|y^n_t|$. It follows from
the convergence of $\{(y^n_t)_{t\in\tT}\}^{+\infty}_{n=1}$ in $\s^2(0,T;\rtn^k)$
that $\sup_{n\in\N}\EXlr{G^2_n(\omega)}<+\infty$. In view of (i) in Lemma
\ref{lem:GnSequenceConvolution} we have that for each $t\in\tT$,
\begin{align*}
  |y^n_t|^2+\intT{t}|z^n_s|^2\dif s
  \leq &\ |\xi|^2-2\intT{t}\langle y^n_s,z^n_s\dif B_s\rangle\\
       &\ +2G_n(\omega)\intT{t}\big(|g(s,0,0)|+kAu(s)(1+|y^n_s|)+kAv(s)(1+|z^n_s|)\big)\dif s.
\end{align*}
It follows from the BDG inequality that
$(\intT[t]{0}\langle y^n_t,z^n_t\dif B_s\rangle)_{t\in\tT}$ is an $(\F_t,\PR)$-martingale.
By the inequalities $2ab\leq a^2+b^2$, $2ab\leq \lambda a^2+b^2/\lambda$
($\lambda:=2k^2A^2\intT{0}v^2(s)\dif s$) and H\"older's inequality we deduce that
for each $n\in\N$,
\begin{align*}
  \EXlr{\intT{0}|z^n_s|^2\dif s}\!
  &\leq\! \EXlr{|\xi|^2}\!+\!\left[1+kA+\lambda+2kA\intT{0}\!\!u(s)\dif s\right]
        \sup_{n\in\N}\!\EXlr{G^2_n(\omega)}
         +\EXlr{\left(\intT{0}\!\!|g(s,0,0)|\dif s\right)^2}\\
  &\hspace{-1cm}+kA\left(\intT{0}\!\!\big(u(s)+v(s)\big)\dif s\right)^2\!\!
   +\frac{k^2A^2}{\lambda}\EXlr{\left(\intT{0}v(s)|z^n_s|\dif s\right)^2}\leq \!C_3+\frac{1}{2}\EXlr{\intT{0}|z^n_s|^2\dif s},
\end{align*}
from which it follows that
\begin{equation}\label{eq:SupNZnNormLeqCLeqInfty}
  \sup_n\EXlr{\intT{0}|z^n_s|^2\dif s}\leq 2C_3<+\infty,
\end{equation}
where $C_3$ is a positive constant and independent of $n$.

On the other hand, by taking expectation in both sides of \eqref{eq:YnAndYmAfterItoFormulaExistencePart},
we have that for each $n$, $m\in\N$,
\begin{equation}\label{eq:ZnZmNormLeqYnYmGnGmWithoutZBds}
  \EXlr{\intT{0}|z^n_s-z^m_s|^2\dif s}
  \leq 2\EXlr{\intT{0}\langle y^n_s-y^m_s,g^n(s,y^n_s,z^n_s)-g^m(s,y^m_s,z^m_s)\rangle\dif s}.
\end{equation}
It follows from (i) in Lemma \ref{lem:GnSequenceConvolution} that
\begin{align*}
  & 2\langle y^n_s-y^m_s,g^n(s,y^n_s,z^n_s)-g^m(s,y^m_s,z^m_s)\rangle\\
  &\leq 4k|y^n_s-y^m_s|[|g(s,0,0)|+Au(s)\big(G_n(\omega)+G_m(\omega)\big)+Av(s)(|z^n_s|+|z^m_s|)+A(u(s)+v(s))].
\end{align*}
Putting the previous inequality into \eqref{eq:ZnZmNormLeqYnYmGnGmWithoutZBds} and
using H\"older's inequality and \eqref{eq:SupNZnNormLeqCLeqInfty} yields that
\begin{align*}
  \EXlr{\intT{0}|z^n_s-z^m_s|^2\dif s}
  \leq &\ 16k\sqrt{\EXlr{\sup_{t\in\tT}|y^n_t-y^m_t|^2}}
        \sqrt{\EXlr{\left(
                      \intT{0}\left[|g(s,0,0)|+A\big(u(s)+v(s)\big)\right]\dif s
                    \right)^2}}\\
       & +32kA\intT{0}u(s)\dif s\sqrt{\sup_{n\in\N}\EXlr{G^2_n(\omega)}}
        \sqrt{\EXlr{\sup_{t\in\tT}|y^n_t-y^m_t|^2}}\\
       & +32kA\sqrt{2C_3\intT{0}v^2(s)\dif s}
         \sqrt{\EXlr{\sup_{t\in\tT}|y^n_t-y^m_t|^2}}.
\end{align*}
Since $\sup_{n\in\N}\EX[G^2_n(\omega)]<+\infty$, $\{(y^n_t)_{t\in\tT}\}_{n=1}^{+\infty}$
converges in $\s^2(0,T;\rtn^k)$ and
\begin{align*}
  & \EXlr{\left(\intT{0}\left[|g(s,0,0)|+A\big(u(s)+v(s)\big)\right]\dif s\right)^2}\\
  & \leq 2\EXlr{\left(\intT{0}|g(s,0,0)|\dif s\right)^2}
  +2\left(\intT{0}A\big(u(s)+v(s)\big)\dif s\right)^2<+\infty,
\end{align*}
we can deduce that,
\begin{equation*}%\label{eq:ZnZmNormConvergeToZero}
  \lim_{n,m\to+\infty}\EXlr{\intT{0}|z^n_s-z^m_s|^2\dif s}=0.
\end{equation*}
Therefore, $\{(z^n_t)_{t\in\tT}\}^{+\infty}_{n=1}$ is a Cauchy sequence in
$\M^2(0,T;\rtn^{k\times d})$. We denote by $(z_t)_{t\in\tT}$ the limit.

\textbf{Step 3.} This step will show that the process $(y_t,z_t)_{t\in\tT}$ is a solution
of BSDE \eqref{eq:BSDEs}.

Now, we have known that for each fixed $t\in\tT$, the sequence $\{y^n_t\}^{+\infty}_{n=1}$ and
$\{\intT{t}z^n_s\dif B_s\}^{+\infty}_{n=1}$ converge in $L^2(\Omega,\F_T,\PR;\rtn^k)$
toward to $y_t$ and $\intT{t}z_s\dif B_s$ respectively. Next let us check the limit
of $g^n(s,y^n_s,z^n_s)$ in BSDE \eqref{eq:BSDEsYnInExistenceProof}. First, for each $t\in\tT$, we have
\begin{align}
  & \EXlr{\left(\intT{t}|g^n(s,y^n_s,z^n_s)-g(s,y_s,z_s)|\dif s\right)^2}\nonumber\\
  &\leq 2\EXlr{\left(\intT{0}|g^n(s,y^n_s,z^n_s)-g(s,y^n_s,z^n_s)|\dif s\right)^2}
    +2\EXlr{\left(\intT{0}|g(s,y^n_s,z^n_s)-g(s,y_s,z_s)|\dif s\right)^2}.
    \label{eq:GnYnZnLeqStepiii}
\end{align}
It follows from (iii) in Lemma \ref{lem:GnSequenceConvolution} that
the first term on the right side of \eqref{eq:GnYnZnLeqStepiii} converges to $0$
as $n\to+\infty$.
Furthermore, by \ref{H:gUniformlyContinuousInY} -- \ref{H:gUniformlyContinuousInZ}
and \eqref{eq:PhiControlAndDivide} we have that for each $n$, $m\in\N$,
\begin{align}
  \EX\!\left[\!\left(\intT{0}\!\!|g(s,y^n_s,z^n_s)\!-\!g(s,y_s,z_s)|\dif s\right)^2\!\right]
  &\!\!\leq\! 2(m+2A)^2\EX\!\left[\!\left(\!
                         \intT{0}\!\!\big(u(s)|y^n_s-y_s|
                         +v(s)|z^n_s-z_s|\big)\dif s
                       \right)^2\!\right]\nonumber\\
  &\quad\!+2\left(\!\intT{0}\!\left[u(s)\rho\left(\frac{2A}{m+2A}\right)
              \!+\!v(s)\phi\left(\frac{2A}{m+2A}\right)\right]\!\dif s
            \right)^2.  \label{eq:GYnZnAndGYZDepartedByYZWithMN}
\end{align}
Note that the second term on the right side of \eqref{eq:GYnZnAndGYZDepartedByYZWithMN}
converges to $0$ as $m\to+\infty$. On the other hand, it follows from H\"older's
inequality that
\begin{align*}
  & \EXlr{\left(\intT{0}\big(u(s)|y^n_s-y_s|+v(s)|z^n_s-z_s|\big)\dif s\right)^2}\\
  & \leq 2\left(\intT{0}u(s)\dif s\right)^2\EXlr{\sup_{t\in\tT}|y^n_t-y_t|^2}
    +2\intT{0}v^2(s)\dif s\EXlr{\intT{0}|z^n_s-z_s|^2\dif s}.
\end{align*}
Thus, by virtue of the fact that $\{(y^n_t,z^n_t)_{t\in\tT}\}_{n=1}^{+\infty}$
is a Cauchy sequence in $\s^2(0,T;\rtn^k)\times \M^2(0,T;\rtn^{k\times d})$,
taking $n\to+\infty$ and then $m\to +\infty$ in \eqref{eq:GYnZnAndGYZDepartedByYZWithMN}
and taking $n\to+\infty$ in \eqref{eq:GnYnZnLeqStepiii} yield that for
each $t\in\tT$,
\begin{equation*}
  \lim_{n\to+\infty}\EXlr{\left|\intT{t}g^n(s,y^n_s,z^n_s)\dif s
  -\intT{t}g(s,y_s,z_s)\dif s\right|^2}=0.
\end{equation*}
Consequently, noticing that $(y_t)_{t\in\tT}$ is a continuous process, by passing to the limit in
BSDE \eqref{eq:BSDEsYnInExistenceProof} we deduce that $\prs$,
\begin{equation*}
  y_t=\xi+\intT{t}g(s,y_s,z_s)\dif s-\intT{t}z_s\dif B_s,\quad t\in\tT,
\end{equation*}
which means that $(y_t,z_t)_{t\in\tT}\in\s^2(0,T;\rtn^k)\times\M^2(0,T;\rtn^{k\times d})$
is a solution of BSDE \eqref{eq:BSDEs}. \hfill$\square$

\begin{proof}[{\bf Proof of Lemma \ref{lem:GnSequenceConvolution}}]
  For each $i=1,\ldots,k$, by \ref{H:gUniformlyContinuousInY} -- \ref{H:ithComponent}
  we deduce that $\pts$, for each $n\in\N$, $y$, $p\in\rtn^k$ and $z$, $q\in\rtn^{k\times d}$,
  \begin{align}
    & g_i(t,p,{}^iq)+(n+A)\big(u(t)|p-y|+v(t)|{}^iq-{}^iz|\big)\nonumber\\
    & \geq g_i(t,0,0)-\big(Au(t)(1+|p|)+Av(t)(1+|{}^iq|)\big)+Au(t)(|p|-|y|)+Av(t)(|{}^iq|-|{}^iz|)\nonumber\\
    %& \geq g_i(t,0,0)-Au(t)(1+|p|)-Av(t)(1+|q|)+A(|u|-|y|)+A(|q|-|z|)\nonumber\\
    & \geq g_i(t,0,0)-Au(t)(1+|y|)-Av(t)(1+|{}^iz|).\label{eq:GIthCompontentGeqGtzero}
  \end{align}
  Thus, for each $n\in\N$, $i=1,\cdots,k$, $y\in\rtn^k$ and $z\in\rtn^{k\times d}$,
  we can define the following $(\F_t)$-progressively measurable function:
  \begin{equation}\label{eq:GIthCompontentNDefineByConvolution}
    g_i^n(t,y,z)=\inf_{(p,q)\in\rtn^k\times\rtn^{k\times d}}
    \left\{g_i(t,p,{}^iq)+(n+A)\big(u(t)|p-y|+v(t)|{}^iq-{}^iz|\big)\right\},
  \end{equation}
  and it depends only on the $i$th row of $z$. Obviously,
  $g^n_i(t,y,{}^iz)\leq g_i(t,y,{}^iz)$, and it follows from
  \eqref{eq:GIthCompontentGeqGtzero} that
  \begin{equation*}
    g^n_i(t,y,{}^iz)\geq -|g_i(t,0,0)|-Au(t)(1+|y|)-Av(t)(1+|{}^iz|).
  \end{equation*}
  Hence, for each $n\in\N$, $g^n_i$ is well defined and (i) holds true with
  setting $g^n:=(g_1^n,g_2^n,\cdots,g_k^n)$.

  Furthermore, it follows from \eqref{eq:GIthCompontentNDefineByConvolution} that
  \begin{equation*}
    g_i^n(t,y,{}^iz)=\inf_{(\overline{p},\overline{q})\in\rtn^k\times\rtn^{k\times d}}
    \left\{g_i(t,y-{\overline p},{}^iz-{}^i{\overline q})+(n+A)\big(u(t)|\overline{p}|+v(t)|{}^i\overline{q}|\big)\right\}.
  \end{equation*}
  Thus, in view of \ref{H:gUniformlyContinuousInY} -- \ref{H:gUniformlyContinuousInZ}
  and the following basic inequality
  \begin{equation}\label{eq:BasicInequalityInfLeqSup}
    \left|\inf_{x\in D} f_1(x)-\inf_{x\in D} f_2(x)\right|\leq
    \sup_{x\in D} |f_1(x)-f_2(x)|,
  \end{equation}
  we have, $\pts$, for each $n\in\N$, $i=1,\cdots,k$, $y_1$, $y_2\in\rtn^k$ and
  $z_1$, $z_2\in\rtn^{k\times d}$,
  \begin{align*}
    |g_i^n(t,y_1,{}^iz_1)-g_i^n(t,y_2,{}^iz_2)|
    & \leq \sup_{(\overline{p},\overline{q})\in\rtn^k\times\rtn^{k\times d}}
    |g_i(t,y_1-{\overline{p}},{}^iz_1-{}^i{\overline{q}})-g_i(t,y_2-{\overline{p}},{}^iz_2-{}^i{\overline{q}})|\\
    & \leq u(t)\rho(|y_1-y_2|)+v(t)\phi(|{}^iz_1-{}^iz_2|),
  \end{align*}
  which means that \ref{H:gUniformlyContinuousInY} -- \ref{H:gUniformlyContinuousInZ}
  hold true for $g^n_i$.

  In the sequel, it follows from \eqref{eq:GIthCompontentNDefineByConvolution} and
  \eqref{eq:BasicInequalityInfLeqSup} that, $\pts$, for each $n\in\N$, $i=1,\cdots,k$,
  $y_1$, $y_2\in\rtn^k$ and $z_1$, $z_2\in\rtn^{k\times d}$,
  \begin{align*}
    & |g_i^n(t,y_1,{}^iz_1)-g_i^n(t,y_2,{}^iz_2)|\\
    & \leq \sup_{(p,q)\in\rtn^k\times\rtn^{k\times d}}
      |(n+A)(u(t)|p-y_1|+v(t)|{}^iq-{}^iz_1|)-(n+A)(u(t)|p-y_2|+v(t)|{}^iq-{}^iz_2|)|\\
    & \leq (n+A)\big(u(t)|y_1-y_2|+v(t)|{}^iz_1-{}^iz_2|\big).
  \end{align*}
  Hence, \ref{A:gNonuniformlyLipschitzContinuousInYAndZ} is right for $g^n_i$.

  Finally, for each $n\in\N$, $i=1,\cdots,k$, $y\in\rtn^k$, $z\in\rtn^{k\times d}$ and $t\in\tT$, let
  \begin{equation*}
    \h_{n,i,t}(y,z):=
    \left\{
      (p,q)\in\rtn^k\times\rtn^{k\times d}:u(t)|p-y|+v(t)|{}^iq-{}^iz|>\frac{2A}{n}\big(u(t)+v(t)\big)
    \right\},
  \end{equation*}
  then
  \begin{equation*}
    \h_{n,i,t}^c(y,z)=
    \left\{
      (p,q)\in\rtn^k\times\rtn^{k\times d}:u(t)|p-y|+v(t)|{}^iq-{}^iz|\leq\frac{2A}{n}\big(u(t)+v(t)\big)
    \right\}.
  \end{equation*}
  For each $n\in\N$, $i=1,\cdots,k$, $y\in\rtn^k$, $z\in\rtn^{k\times d}$, $t\in\tT$
  and $(p,q)\in\h_{n,i,t}(y,z)$,
  it follows from \ref{H:gUniformlyContinuousInY} -- \ref{H:gUniformlyContinuousInZ}
  that
  \begin{align*}
    g_i(t,p,{}^iq)\!+\!(n+A)\big(u(t)|p-y|\!+\!v(t)|{}^iq-{}^iz|\big)
    & \geq g_i(t,y,{}^iz)\!+\!n\big(u(t)|p-y|\!+\!v(t)|{}^iq-{}^iz|\big)\!-\!A\big(u(t)\!+\!v(t)\big)\\
    & >g_i(t,y,{}^iz)+A\big(u(t)+v(t)\big),\quad \prs.
  \end{align*}
  Then, since $g^n_i(t,y,{}^iz)\leq g_i(t,y,{}^iz)$, we have that $\pts$, for each $n\in\N$,
  $i=1,\cdots,k$, $y\in\rtn^k$ and $z\in\rtn^{k\times d}$,
  \begin{equation}\label{eq:GnIthCompontentDefineByHc}
    g_i^n(t,y,{}^iz)=\inf_{(p,q)\in\h_{n,i,t}^c(y,z)}
    \left\{g_i(t,p,{}^iq)+(n+A)\big(u(t)|p-y|+v(t)|{}^iq-{}^iz|\big)\right\}.
  \end{equation}
  In the sequel, \ref{H:gUniformlyContinuousInY} -- \ref{H:ithComponent}
  and \eqref{eq:GnIthCompontentDefineByHc} yield that, $\pts$, for each $n\in\N$,
  $i=1,\cdots,k$, $y\in\rtn^k$ and $z\in\rtn^{k\times d}$,
  \begin{align*}
    g_i^n(t,y,{}^iz)&\geq\inf_{(p,q)\in\h_{n,i,t}^c(y,z)}
    \{g_i(t,y,{}^iz)-u(t)\rho(|p-y|)-v(t)\phi(|{}^iq-{}^iz|)\}
    \geq g_i(t,y,{}^iz)-b_n(t),
  \end{align*}
  where
  \begin{equation*}
    b_n(t)=u(t)\rho\left(\frac{2A}{n}\cdot\frac{u(t)+v(t)}{u(t)}\right)
           +v(t)\phi\left(\frac{2A}{n}\cdot\frac{u(t)+v(t)}{v(t)}\right).
  \end{equation*}
  Thus, $\pts$, for each $n\in\N$, $i=1,\cdots,k$, $y\in\rtn^k$ and
  $z\in\rtn^{k\times d}$,
  \begin{equation*}
    0\leq g_i(t,y,{}^iz)-g_i^n(t,y,{}^iz)\leq b_n(t).
  \end{equation*}
  It is clear that $b_n(t)\downarrow0$ as $n\to +\infty$ for each $t\in\tT$.
  Since $\rho(\cdot)$ and $\phi(\cdot)$ are at most linear-growth, we have
  \begin{equation*}
    b_n(t)\leq Au(t)+\frac{2A^2}{n}\big(u(t)+v(t)\big)+Av(t)
               +\frac{2A^2}{n}\big(u(t)+v(t)\big)
          \leq (A+4A^2)\big(u(t)+v(t)\big),
  \end{equation*}
  from which and the Lebesgue dominated convergence theorem it follows that
  $\intT{0}b_n(t)\dif t\to 0$ as $n\to +\infty$.

  Thus, the sequence $g^n:=(g_1^n,g_2^n,\cdots,g_k^n)$ is just the one we desire.
  The proof is complete.
\end{proof}

\begin{rmk}
  Note that if there exists a constant $K>0$ such that the function $\phi(\cdot)$
  appearing in \ref{H:gUniformlyContinuousInZ} satisfies $\phi(x)\leq Kx$ for all
  $x\in\rtn^+$, then the condition $\intT{0}\big(v(t)+v^2(t)\big)\dif t<+\infty$
  in Theorem \ref{thm:MainResult} can be weakened to $\intT{0}v^2(t)\dif t<+\infty$
  as in Lemma \ref{lem:ResultInChenWang}.
\end{rmk}

%\appendix
\section*{Appendix: Some proofs of results on DBDEs}
For the convenience of readers, we would like to provide here the proofs of those results
on deterministic backward stochastic differential equations introduced
in Section \ref{sec:Preliminaries}.

\begin{proof}[{\bf Proof of Proposition \ref{pro:DBDEResult}}]
  For each $\beta(\cdot):\rtn^+\mapsto\rtn^+$ such that $\intT{0}\beta(t)\dif t<+\infty$,
  let $\CH_{\beta(\cdot)}$ denote the set of the continuous functions
  $(y_t)_{t\in\tT}$ such that
  $$\|y\|_{\beta(\cdot)}:=
    \left(
      \sup_{t\in\tT}\left[\me^{-\intT{t}\beta(r)\dif r}|y_t|^2\right]
    \right)^{1/2}<+\infty.$$
  It is easy to verify that $\CH_{\beta(\cdot)}$ is a Banach space. Note that for any
  $y_t\in\CH_{\beta(\cdot)}$, in view of
  \ref{B:fNonuniformlyLipschitzContinuous} and \ref{B:fOnlyIntegrable},
  \begin{align*}
    \intT{0}|f(t,y_t)|\dif t
    & \leq \intT{0}u(t)|y_t|\dif t+\intT{0}|f(t,0)|\dif t\\
    & \leq \intT{0}u(t)\dif t\cdot
           \left[
             \me^{\intT{0}\beta(r)\dif r}\|y\|_{\beta(\cdot)}
           \right]
           +\intT{0}|f(t,0)|\dif t<+\infty.
  \end{align*}
  For any $y_t\in\CH_{\beta(\cdot)}$, define
  \begin{equation*}%\label{eq:YWithDBDEToEstablishStrictContraction}
    Y_t:=\delta+\intT{t}f(s,y_s)\dif s,\quad t\in\tT.
  \end{equation*}
  Then we have
  \begin{align*}
    \|Y\|_{\beta(\cdot)}
    =\sqrt{\sup_{t\in\tT}
             \left[
               \me^{-\intT{t}\beta(r)\dif r}|Y_t|^2
             \right]} \leq \sup_{t\in\tT}|Y_t|
    \leq |\delta|+\intT{0}|f(t,y_t)|\dif t<+\infty.
  \end{align*}
  Thus, we have constructed a mapping $\Phi:\CH_{\beta(\cdot)}\mapsto\CH_{\beta(\cdot)}$
  such that $\Phi(y_t)=Y_t$. Next we prove that this mapping is strictly contractive
  when $\beta(\cdot)$ is chosen appropriately.

  Take $y^1_t$, $y^2_t\in\CH_{\beta(\cdot)}$ and assume that $\Phi(y^1_t)=Y^1_t$,
  $\Phi(y^2_t)=Y^2_t$. Let us set $\hat Y_t:=Y^1_t-Y^2_t$, $\hat y_t:=y^1_t-y^2_t$.
  Then we have
  \begin{align*}
    \dif\left[\me^{-\intT{s}\beta(r)\dif r}|\hat Y_s|^2\right]
    & =\me^{-\intT{s}\beta(r)\dif r}
       \big(\beta(s)|\hat Y_s|^2\dif s+2\hat Y_s\dif \hat Y_s\big)\\
    & =\me^{-\intT{s}\beta(r)\dif r}
       \big(\beta(s)|\hat Y_s|^2\dif s-2\hat Y_s(f(s,y^1_s)-f(s,y^2_s))\dif s\big),
  \end{align*}
  from which it follows that
  \begin{equation*}
    \me^{-\intT{t}\beta(r)\dif r}|\hat Y_t|^2
    =\intT{t}\me^{-\intT{s}\beta(r)\dif r}
     \left[
       2\hat Y_s
       \big(
         f(t,y^1_s)-f(t,y^2_s)
       \big)
       -\beta(s)|\hat Y_s|^2
     \right]\dif s.
  \end{equation*}
  By \ref{B:fNonuniformlyLipschitzContinuous} and the inequality
  $2ab\leq\lambda a^2+b^2/\lambda$ $(\lambda>0)$ we get that
  \begin{align*}
    \me^{-\intT{t}\beta(r)\dif r}|\hat Y_t|^2
    & \leq \intT{t}\me^{-\intT{s}\beta(r)\dif r}
           \left[
             2\sqrt{u(s)}|\hat Y_s|\sqrt{u(s)}|\hat y_s|-\beta(s)|\hat Y_s|^2
           \right]\dif s \nonumber\\
    & \leq \intT{t}\me^{-\intT{s}\beta(r)\dif r}
           \left[
             \big(
               \lambda u(s)-\beta(s)
             \big) |\hat Y_s|^2+\frac{u(s)}{\lambda}|\hat y_s|^2
           \right]\dif s, \label{eq:eBetaHatYLeqLambda}
  \end{align*}
  from which it follows that, with choosing $\lambda>4\intT{0}u(s)\dif s$ and
  $\beta(s)=\lambda u(s)$,
  \begin{align*}
    \|\hat Y\|^2_{\beta(\cdot)}
    \leq \frac{1}{\lambda}\intT{0}\me^{-\intT{s}\beta(r)\dif r}u(s)|\hat y_s|^2\dif s
    \leq \frac{1}{\lambda}\sup_{t\in\tT}
           \left[
             \me^{-\intT{t}\beta(r)\dif r}|\hat y_t|^2
           \right]\intT{0}u(s)\dif s
      <\frac{1}{4}\|\hat y\|^2_{\beta(\cdot)}.
  \end{align*}
  So $\|\hat Y\|_{\beta(\cdot)}<\frac{1}{2}\|\hat y\|_{\beta(\cdot)}$,
  which implies that $\Phi$ is a contractive mapping from $\CH_{\beta(\cdot)}$
  to $\CH_{\beta(\cdot)}$. Then the conclusion follows from the fixed point theorem
  immediately.
\end{proof}

Proposition \ref{pro:DBDESolutionsDefinedRecursively} follows directly from the
proof Proposition \ref{pro:DBDEResult}.
%From the proof of Theorem \ref{pro:DBDEResult}, the following Proposition
%\ref{pro:DBDESolutionsDefinedRecursively} is immediate.

\begin{proof}[{\bf Proof of Proposition \ref{pro:DBDESolutionComparison}}]
  For each $t\in\tT$, let us set
  \begin{equation*}
    a(t):=
    \begin{cases}
      \displaystyle \frac{f(t,y_t)-f(t,y'_t)}{y_t-y'_t}, & y_t\neq y'_t;\\
      \displaystyle 0, & y_t=y'_t,
    \end{cases}
  \end{equation*}
  $b(t):=f(t,y'_t)-f'(t,y'_t)$ and $\hat y_t:=y_t-y'_t$. From
  \ref{B:fNonuniformlyLipschitzContinuous} and \ref{B:fOnlyIntegrable} we can deduce that
  $\intT{0}|a(s)|\dif s\leq\intT{0}u(s)\dif s<+\infty$ and $\intT{0}b(s)\dif s<+\infty$.
  Then we have
  \begin{align*}
    \hat y_t
    %& =\delta-\delta'+\intT{t}\big(f(s,y_s)-f'(s,y'_s)\big)\dif s\\
    & =\delta-\delta'+\intT{t}\big(f(s,y_s)-f(s,y'_s)+f(s,y'_s)-f'(s,y'_s)\big)\dif s\\
    & =\delta-\delta'+\intT{t}a(s)\hat y_s\dif s+\intT{t}b(s)\dif s,\quad t\in\tT.
  \end{align*}
  Consequently, in view of the conditions that $\delta\geq\delta'$ and $b(t)\geq 0$ for
  each $t\in\tT$,
  \begin{equation*}
    \hat y_t=\me^{\intT{t}a(s)\dif s}
             \left[
               \delta-\delta'+\intT{t}b(s)\me^{\intT{s}a(r)\dif r}\dif s
             \right]\geq 0, \quad t\in\tT.
  \end{equation*}
  Furthermore, if $\delta>\delta'$, then $\hat y_t>0$ for each $t\in\tT$.
  The proof is then completed.
\end{proof}

\begin{proof}[{\bf Proof of Proposition \ref{pro:DBDEsWithTwoFurthermoreResults}}]
  Let $\varphi_n(x):=\inf_{y\in\rtn}\{\overline{\varphi}(y)+n|x-y|\}$ with
  $\overline{\varphi}(y)=\varphi(|y|)$ for $y\in\rtn$, then it follows from Lemma
  \ref{lem:LepeltierMartinLemma} and Proposition \ref{pro:DBDEResult} that $\varphi_n(x)$
  is well defined on $\rtn$ for each $n\geq a$, and for each $n\geq a$, the
  following two DBDEs
  \begin{equation}\label{eq:BDEfnDeltaPhin}
    y^{n,\delta}_t=\delta+\intT{t}u(s)\varphi_n(y^{n,\delta}_s)\dif s,\quad t\in\tT,
  \end{equation}
  and
  \begin{equation*}
    \overline{y}^\delta_t=\delta+\intT{t}\big(au(s)\overline{y}^\delta_s+bu(s)\big)\dif s,\quad t\in\tT,
  \end{equation*}
  have, respectively, unique solutions $y^{n,\delta}_t$ and $\overline{y}^\delta_t$
  with $\sup_{t\in\tT}|y^{n,\delta}_t|<+\infty$ and
  $\sup_{t\in\tT}|\overline{y}^\delta_t|<+\infty$. Clearly, $y^{n,\delta}_t\geq 0$
  and $\overline{y}^\delta_t\geq 0$ for each $t\in\tT$ and $n\geq a$.
  By Proposition \ref{pro:DBDESolutionComparison} and the fact that
  $\varphi_n\leq\varphi_{n+1}$, we have $y^{n,\delta}\leq y^{n+1,\delta}\leq\overline{y}^\delta$.
  Therefore, for each $t\in\tT$, the limit of the sequence $\{y^{n,\delta}_t\}_{n=a}^{+\infty}$
  must exist, we denote it by $y^\delta_t$. In view of (i) and (iv) in Lemma
  \ref{lem:LepeltierMartinLemma}, using Lebesgue's dominated convergence theorem
  we can obtain that
  \begin{equation*}
    \lim_{n\to+\infty}\intT{0}u(s)\varphi_n(y^{n,\delta}_s)\dif s
    =\intT{0}u(s)\overline{\varphi}(y^\delta_s)\dif s
    =\intT{0}u(s)\varphi(y^\delta_s)\dif s.
  \end{equation*}
  Thus, by passing to the limit in both sides of DBDE \eqref{eq:BDEfnDeltaPhin}
  we deduce that
  \begin{equation*}
    y^\delta_t=\delta+\intT{t}u(s)\varphi(y^\delta_s)\dif s, \quad t\in\tT,
  \end{equation*}
  which means that $y^\delta_t$ is a solution of DBDE \eqref{eq:DBDEFurthermoreResults}.

  Let us now suppose that $\delta>0$ and $\varphi(x)>0$ for all $x>0$. For each $z\geq\delta$,
  set $G(z):=\intT[1]{z}\varphi^{-1}(x)\dif x$. It is easy to see that
  $-\infty=G(+\infty)<G(z_1)<G(z_2)<G(\delta)$ for each $z_1> z_2>\delta$. Then
  the inverse function of $G(z)$ must exist, we denote it by $G^{-1}(u)$ for
  $u\leq G(\delta)$. Let $y^\delta_t$ be a solution of
  DBDE \eqref{eq:DBDEFurthermoreResults}.
  It is obvious that $y^\delta_t\geq\delta$ and $\dif G(y^\delta_t)=u(t)\dif t$.
  Hence, $G(y^\delta_T)-G(y^\delta_t)=\intT{t}u(s)\dif s$, which implies that
  $y^\delta_t=G^{-1}\big(G(\delta)-\intT{t}u(s)\dif s\big)$ for each $t\in\tT$.
  The proof of (i) is then completed.

  Finally, we prove that (ii) is also right. Assume that $\delta=0$ and $\varphi\in\BS$
  with $\int_{0^+}\varphi^{-1}(x)\dif x=+\infty$. For each $z>0$, set
  $H(z):=\intT[1]{z}\varphi^{-1}(x)\dif x$. It is clear that
  $-\infty=H(+\infty)<H(z_1)<H(z_2)<H(0)=+\infty$ for each $z_1>z_2>0$. Then the
  inverse function of $H(z)$ must exist, we denote it by $H^{-1}(u)$ for each $u\in\rtn$.
  Now let $y^0_t$ be a continuous solution of DBDE \eqref{eq:DBDEFurthermoreResults}.
  Then $y^0_t\geq 0$ and $\dif H(y^0_t)=u(t)\dif t$. Hence, for each $0\leq t\leq t_1<T$,
  \begin{equation*}
    y^0_t=H^{-1}\left(H(y^0_{t_1})-\intT[t_1]{t}u(s)\dif s\right).\smallskip
  \end{equation*}
  Furthermore, noticing that $H(y^0_{t_1})\to H(y^0_T)=H(0)=+\infty$ as $t_1\to T$ and
  $H^{-1}(+\infty)=0$, we know that $y^0_t=0$ for each $t\leq T$. Consequently,
  DBDE \eqref{eq:DBDEFurthermoreResults} has a unique solution $y_t\equiv 0$.
\end{proof}

\bibliographystyle{elsarticle-harv}
\bibliography{bsdes}

%\newpage
%\renewcommand{\baselinestretch}{1}
%\pagenumbering{arabic}

%\section*{Some proofs of results about DBDEs}\label{App:ProofsOfResultsAboutDBDEs}

\end{document}